\documentstyle[aps,epsfig]{revtex}
\begin{document}
\twocolumn
\title{Hybrid-Cubic-Rational Semi-Lagrangian Method with the Optimal Mixing}
\author{Masato Ida}
\address{Satellite Venture Business Laboratory, Gunma University, 1--5--1 Tenjin-cho, Kiryu-shi, Gunma 376-8515, Japan\\
E-mail : ida@vbl.gunma-u.ac.jp\\ TEL : +81-277-30-1126   FAX : +81-277-30-1121}
\maketitle
\begin{abstract}
A semi-Lagrangian method for advection equation with hybrid cubic-rational 
interpolation is introduced. In the present method, the spatial profile of 
physical quantities is interpolated with a combination of a cubic and a 
rational function. For achieving both high accuracy and convexity preserving 
of solution, the two functions are mixed in the optimal ratio which is given 
theoretically. Accuracy and validity of this method is demonstrated with 
some numerical experiments.\\ \\
\noindent
{\bf Key Words:} Numerical method, Advection, Semi-Lagrangian method, 
Interpolation, Cubic function, Rational function, Convexity preserving.
\end{abstract}

\section{INTRODUCTION}
Basically, the CIP method \cite{ref1,ref2} is a numerical method for an 
advection equation,
\[
\frac{{\partial f}}{{\partial t}} + u\frac{{\partial f}}{{\partial x}} = 0
\]
whose solution is expressed as
\begin{equation}
\label{eq1}
f(x,t + \Delta t) = f(X(x,t),t),
\end{equation}
where $X$ is the trajectory of fluid particle, which is located at 
$x$ at the time $t + \Delta t$,
\begin{equation}
\label{eq2}
X(x,t) = x + \int_{t}^{t - \Delta t} {u(X(x,\tau ),\tau )d\tau }.
\end{equation}
In a semi-Lagrangian scheme like the CIP, the solution (\ref{eq1}) is solved as an 
interpolation problem \cite{ref3,ref1}. In the CIP scheme, an Hermite cubic 
expansion function is used to interpolate $f$ at the time $t$. The 
quantity and its first spatial derivative defined at each grid points are 
updated as to obey, respectively, eq.~(\ref{eq1}) and its spatial derivative, i.e.,
\[
\frac{{\partial f(x,t + \Delta t)}}{{\partial x}} = \frac{{\partial X(x,t)}}
{{\partial x}}\frac{{\partial f(X(x,t),t)}}{{\partial X}}
\]
Generally, the integration in eq.~(\ref{eq2}) is solved by assuming that $u$ 
is locally constant as
\begin{eqnarray}
X(x,t) &\approx& x + u(x,t)\int_{t}^{t - \Delta t} {d\tau } \nonumber\\
&=& x - u(x,t)\Delta t.\nonumber
\end{eqnarray}
With this, the solutions (\ref{eq1}) and (\ref{eq2}) are expressed, 
respectively, as
\begin{equation}
\label{eq3}
f(x,t + \Delta t) = f(x - u(x,t)\Delta t,t),
\end{equation}
\begin{equation}
\label{eq4}
\frac{{\partial f(x,t + \Delta t)}}{{\partial x}} = [1 - \frac{{\partial 
u(x,t)}}{{\partial x}}\Delta t]\frac{{\partial f(x - u(x,t)\Delta 
t,t)}}{{\partial X}}.
\end{equation}

In the last decade, various kinds of extension and improvement have been 
adopted to this method. In 1990, Yabe et al extended this to multidimensions 
without time-splitting technique by employing a multidimensional cubic 
expansion function \cite{ref4}. In 1991, Kondoh extended the 1D CIP to a 
5th-order advection method and, furthermore, proposed a solver for parabolic 
equations by extending the basic concept of the CIP \cite{ref5}. In the same 
year, Aoki and Yabe improved the multidimensional one by modifying the 
multidimensional expansion function and proposed two alternative formulae 
\cite{ref6,ref7}. In 1994, 
Kondoh extended his approach to a multidimensional parabolic equation and 
general hyperbolic equations \cite{ref8}. In 1995, Ida and Yabe proposed an 
implicit version of the CIP \cite{ref9}. This method is CFL free and can be 
solved directly 
with a marching procedure although it is a 3rd-order method. In the same 
year, Utsumi extended the CIP to a solver for the Euler equations of fluid 
flow without finite-difference technique by employing differential-algebraic 
and Lagrangian-like concepts \cite{ref10}. In 1996, Xiao et al proposed a 
convexity preserving method for the advection equation by replacing the cubic 
function in the CIP scheme with a cubic-rational function \cite{ref11,ref12}. 
In the same year, Ida proposed a high accurate solver for free-surface flow 
problem by coupling the CIP with newly proposed extrapolation scheme 
\cite{ref13,ref14,ref15}. With this 
method, the density discontinuity at material interface is solved without 
any numerical dissipation across the interface. In 1997, by extending the 
Kondoh's approach \cite{ref5}, Aoki proposed high accurate solver for wave 
equation, full Euler equations and others \cite{ref16}. In 1999, Tanaka et al 
proposed an exactly conservative solver for the continuity equation in 
non-conservative form by additionally using the local mass of fluid as a 
dependent variable \cite{ref17,ref18}.

In this paper we discuss on the rational method proposed by Xiao et al. As 
proved theoretically in ref.~11, the rational interpolation method 
suppresses the numerical oscillation which tends to appear in high-order 
solution. However, this method sometime provides more diffusive result than 
that with the classical cubic interpolation. We try to improve its accuracy, 
without any loss of the convexity preserving property, by mixing with the 
cubic interpolation function in the optimal ratio. In Sec.~2, the 
conventional rational method is briefly reviewed and, in Sec.~3, the optimal 
mixing ratio is given theoretically. In Sec.~4, results of some numerical 
experiments are shown for demonstrating the accuracy and the validity of the 
present method.

\section{CONVEXITY PRESERVING SEMI-LAGRANGIAN METHOD}

In the rational method \cite{ref11}, the following cubic-rational interpolation 
function is used:
\begin{equation}
\label{eq5}
CR(\xi ,\gamma) = \frac{{f_{i} + A1_{i} \xi + A2_{i} \xi ^{2} + A3_{i} \xi 
^{3}}}{{1 + \gamma B_{i} \xi }},
\end{equation}
where
\begin{eqnarray}
A1_{i} &=& d_i  + f_i \gamma B_i,\nonumber \\
A2_{i} &=& S_{i} \gamma B_{i} + (S_{i} - d_{i} )/h - A3_{i} h,\nonumber \\
A3_{i} &=& [d_{i} - S_{i} + (d_{i + 1} - S_{i} )(1 + \gamma B_{i} h)]/h^{2},
\nonumber \\
B_{i} &=& \left({\frac{{S_{i} - d_{i}}}{{d_{i + 1} - S_{i}}} - 1} \right)/h,
\nonumber \\
S_{i} &=& (f_{i + 1} - f_{i} )/h,\nonumber \\
\xi &=&  - u_i \Delta t,\nonumber
\end{eqnarray}
$f_i$ and $d_i$ are a physical quantity and its first 
spatial derivative, respectively, ${u_i}$ is the particle velocity 
at $x_i$ assumed as negative here, $\Delta t$ is the time 
interval, $h$ is the grid width assumed as uniform in this paper for 
simplicity and $\gamma $ is a parameter for switching the form of 
the interpolation function. In the case of $\gamma = 1$, eq.~(\ref{eq5}) 
is reduced to a rational formula of
\begin{equation}
\label{eq6}
CR(\xi ,1) = \frac{{f_{i} + R1_{i} \xi + R2_{i} \xi ^{2}}}{{1 + B_{i} \xi 
}},
\end{equation}
where
\begin{eqnarray}
R1_{i} &=& d_{i} + f_{i} B_{i},\nonumber \\
R2_{i} &=& S_{i} B_{i} + (S_{i} - d_{i} )/h.\nonumber
\end{eqnarray}
Unlike rational functions used in data interpolation technique (See 
Ref.~\cite{ref19} for example), the above formula is constructed not only with 
the quantity but also with its derivative. In the case of $\gamma = 0$, on the 
contrary, eq.~(\ref{eq5}) is reduced to a cubic formula of
\begin{equation}
\label{eq7}
CR(\xi ,0) = f_i  + C1_i \xi  + C2_i \xi ^2  + C3_i \xi ^3,
\end{equation}
where
\begin{eqnarray}
C1_{i} &=& d_{i},\nonumber \\
C2_{i} &=& - (2d_{i} + d_{i + 1} - 3S_{i} )/h,\nonumber \\
C3_{i} &=& (d_{i} + d_{i + 1} - 2S_{i} )/h^{2}.\nonumber
\end{eqnarray}
Those rational and cubic formulae satisfy a continuity condition at 
$x_i$ and $x_{i+1}$ expressed as
\[
\left\{ {\begin{array}{*{20}c}
   {CR(0,1) = f_i ,} \hfill  \\
   {CR(h,1) = f_{i + 1} ,} \hfill  \\
   {\partial _x CR(0,1) = d_i ,} \hfill  \\
   {\partial _x CR(h,1) = d_{i + 1} .} \hfill  \\
\end{array}} \right.
\]
The rational formula is adapted in a cell where the data is convex or 
concave, i.e.,
\[
d_i  > S_i  > d_{i + 1}
\]
or
\[
d_{i} < S_{i} < d_{i + 1}.
\]
This rational interpolation function preserves the convexity of solution. 
For the other data, interpolation function is switched to the formula of 
eq.~(\ref{eq7}) which corresponds to the conventional one of CIP \cite{ref1} 
and provides purely 3rd-order solution.

With those interpolation functions, in this paper, we propose a hybrid 
method of making use only of their superior characteristics by mixing them 
optimally. In ref.~12, Xiao et al proposed an additional switching technique 
shown as
\begin{equation}
\label{eq8}
\gamma  = \left\{ {\begin{array}{*{20}c}
   {1,\quad} \hfill & {{\rm for}\;d_i  \cdot d_{i + 1}  < 0,} \hfill  \\
   {0,\quad} \hfill & {\rm otherwise.} \hfill  \\
\end{array}} \right.
\end{equation}
This means that the rational function is applied only in the cell which 
includes a turning point of the gradient. While this procedure modifies the 
dissipation property of the rational method, this breaks the preserving of 
convexity as shown in Sec.~4. The optimal mixing technique being proposed in 
this paper would achieve improvement of accuracy without any loss of the 
convexity preserving property.

For the convenience of the following discussion, we rearrange eqs.~(\ref{eq6}) 
and (\ref{eq7}) as
\begin{equation}
\label{eq9}
R(k) = f_{i} + d_{i} hk + \frac{{P_{i} ^{2}k^{2}}}{{Q_{i} + (P_{i} - Q_{i})k}},
\end{equation}
and
\begin{equation}
\label{eq10}
C(k) = f_{i} + d_{i} hk + (2P_{i} - Q_{i} )k^{2} + (Q_{i} - P_{i} )k^{3},
\end{equation}
respectively, where
\begin{eqnarray}
P_i &=& (S_i  - d_i )h,\nonumber \\
Q_{i} &=& (d_{i + 1} - S_{i} )h,\nonumber
\end{eqnarray}
and
\[
k \equiv \xi /h = - u_{i} \Delta t/h
\]
is the local Courant number and
\[
k \in [0,1]
\]
because of the CFL condition.

\section{HYBRID CUBIC-RATIONAL METHOD WITH THE OPTIMAL MIXING}

\subsection{The optimal mixing of the two interpolation functions under the 
convexity-preserving condition}

We start from a combination of the rational and the cubic functions shown as
\begin{equation}
\label{eq11}
F(k) = \alpha R(k) + (1 - \alpha )C(k),
\end{equation}
where $\alpha $ is a weighting parameter whose range is limited as 
$\alpha \in [0,1]$. For large $\alpha $, i.e., $\alpha \sim 1$, 
less oscillatory solution may be 
expected because the rational function becomes dominant and, for small 
$\alpha $, i.e., $\alpha \sim 0$, high-order solution may be 
expected because the cubic one becomes dominant. By determining 
$\alpha $ properly, convexity-preserving high-accurate method may 
be produced. We discuss below how to determine the weighting parameter.

While numerous proofs on some characteristics of the rational interpolation 
method have been done in Ref.~\cite{ref11}, the roof of all of them is a 
property of the rational function. The property is expressed as that, for the 
convex data, the interpolation function is convex between the given interval 
and, for the concave data, it is concave. Namely,
\[
\frac{{\partial^{2}R(k)}}{{\partial x^{2}}} \le 0,\quad {\rm for}\;k \in [1,0]
\]
is true for the convex data of
\[
P_{i} < 0\;{\rm and}\;Q_{i} < 0,
\]
and
\[
\frac{{\partial ^{2}R(k)}}{{\partial x^{2}}} \ge 0,\quad{\rm for}\;k \in [1,0]
\]
is true for the concave data of
\[
P_{i} > 0\; {\rm and} \; Q_i > 0.
\]
For achieving both convexity preserving and better resolution, we make the 
weighting parameter the minimum value which satisfies the above condition. 
The following discussion will be limited to the case where the data is 
convex or concave since the rational interpolation is adopted only in the 
case. For the other data, the cubic interpolation is adopted.

The second spatial derivative of eq.~(\ref{eq11}) is shown as
\begin{eqnarray}
\label{eq12}
\frac{\partial ^2 F(k)}{\partial x^2 }&& = \alpha \frac{2P_i ^2 Q_i ^2}
{h^2 [(P_i - Q_i )k + Q_i]^3} \nonumber \\
&&+ (1 - \alpha )\frac{2}{h^2} [3(Q_i - P_i)k  + 2 P_i - Q_i].
\end{eqnarray}
By introducing $H_{i}$ defined as
\[
H_{i} \equiv \frac{Q_{i}}{P_{i}},
\]
eq.~(\ref{eq12}) is rewritten as
\begin{eqnarray}
\frac{{\partial ^2 F(k)}}{{\partial x^2 }} = \frac{2}{h^2}Q_i
\left\{ \alpha \frac{{H_i }}{{[(1 - H_i )k + H_i ]^3 }} \right. \nonumber \\
\left. + (1 - \alpha ) \frac{{3(H_i - 1)k + 2 - H_i }}{{H_i}} \right\}.
\nonumber
\end{eqnarray}
For preserving convexity of solution, as mentioned above, the following 
condition must be satisfied:
\[
\frac{{\partial ^{2}F(k)}}{{\partial x^{2}}} \le 0,\quad {\rm for}\;P_{i} < 0
\;{\rm and}\;Q_{i} < 0,
\]
or
\[
\frac{{\partial ^2 F(k)}}{{\partial x^2 }} \ge 0,\quad {\rm for}\;P_i > 0\;
{\rm and}\;Q_i > 0.
\]
This condition is expressed by an inequality of
\begin{eqnarray}
\label{eq13}
{{\frac{{\partial ^2 F(k)}}{{\partial x^2 }}} \mathord{\left/
 {\vphantom {{\frac{\partial ^2 F(k)}{\partial x^2 }}
{\left( {\frac{{h^2 }}{2}Q_i } \right)}}} \right.
\kern-\nulldelimiterspace} {\left( {\frac{2}{h^2}Q_i } \right)}}
= \alpha \frac{H_i}{[(1 - H_i )k + H_i ]^3 } \nonumber \\
+ (1 - \alpha )\frac{{3(H_i - 1)k + 2 - H_i }}{{H_i }} \ge 0,
\end{eqnarray}
because the sigh of ${{\partial ^{2}F(k)} \mathord{\left/ {\vphantom 
{{\partial ^{2}F(k)} {\partial x^{2}}}} \right. \kern-\nulldelimiterspace} 
{\partial x^{2}}}$ and $Q_i$ should be the same. Furthermore, the 
inequality (\ref{eq13}) can be rewritten as
\begin{equation}
\label{eq14}
\alpha H_{i} ^{2} + (1 - \alpha )G_{i} ^{3}[3(H_{i} - 1)k + 2 - H_{i}] \ge 0,
\end{equation}
where
\[
G_{i} \equiv (1 - H_{i} )k + H_{i},
\]
because of
\begin{equation}
\label{eq15}
H_i  > 0
\end{equation}
(Remind that the signs of $P_i$ and $Q_i$ are the same) and
\begin{equation}
\label{eq16}
G_{i} \ge {\rm min}(1,H_{i} ) > 0, \quad {\rm for}\;k \in [0,1].
\end{equation}
Here we introduce $K_i(k)$ which is defined as
\[
K_{i} (k) \equiv G_{i} ^{3}[3(H_{i} - 1)k + 2 - H_{i} ].
\]
With this function, eq.~(\ref{eq14}) is rewritten as
\begin{equation}
\label{eq17}
\alpha H_{i} ^{2} + (1 - \alpha )K_{i} (k) \ge 0.
\end{equation}
From the inequality (\ref{eq16}) and
\begin{eqnarray}
3(H_{i} - 1)k + 2 - H_{i} &\ge& {\rm min}(2 - H_{i} ,2H_{i} - 1),\nonumber \\
{\rm for}\;k &\in& [0,1],\nonumber
\end{eqnarray}
it is known that the inequality (\ref{eq17}) is always true for
\[
1/2 \le H_i \le 2
\]
because of $K_{i} (k) \ge 0$. In this case, we can use $\alpha = 0$, i.e., 
the purely cubic interpolation function which provides high-order solution.

For a while, the following discussion will be limited to the case of
\[
H_{i} > 2.
\]
The first and the second spatial derivatives of $K_{i}(k)$ 
respect to $k$ are
\begin{equation}
\label{eq18}
\frac{{\partial K_{i} (k)}}{{\partial k}} = 6(1 - H_{i})^{2}G_{i} ^{2}(1 - k)
\end{equation}
and
\begin{equation}
\label{eq19}
\frac{{\partial ^{2}K_{i} (k)}}{{\partial k^{2}}} = 6(1 - H_{i})^{2}G_{i} 
(2 - 3G_{i} ),
\end{equation}
respectively. From eq.~(\ref{eq19}) with an inequality of
\[
G_i \ge \min (1,H_i) = 1,\quad {\rm for}\;H_{i} \ge 2\;{\rm and}\;k \in [0,1],
\]
we know that
\[
\frac{\partial ^{2}K_{i}(k)}{\partial k^{2}} < 0,\quad {\rm for}\;k \in [0,1],
\]
which means that ${{\partial K_{i} (k)} \mathord{\left/ {\vphantom 
{{\partial K_{i} (k)} {\partial k}}} \right. \kern-\nulldelimiterspace} 
{\partial k}}$ is monotonicaly decreasing for $k \in [0,1]$. Furthermore, from 
eq.~(\ref{eq18}), we know that
\[
\frac{{\partial K_{i} (0)}}{{\partial k}} = 6(1 - H_{i} )^{2}H_{i} ^{2} > 0
\]
and
\[
\frac{{\partial K_{i} (1)}}{{\partial k}} = 0.
\]
Namely, ${{\partial K_{i} (k)} \mathord{\left/ {\vphantom {{\partial K_{i} 
(k)} {\partial k}}} \right. \kern-\nulldelimiterspace} {\partial k}}$ has 
non-negative value for $k \in [0,1]$.

Those results show that $K_i (k)$ is monotonicaly increasing for $k \in [0,1]$,
and the minimum value of it is
\[
{\rm min}(K_{i} (k)) = K_{i} (0) = H_{i} ^{3}(2 - H_{i} ) < 0.
\]
Thus, the inequality (\ref{eq17}) is reduced to
\[
\alpha H_{i} ^{2} + (1 - \alpha )H_{i} ^{3}(2 - H_{i} ) \ge 0.
\]
From this, we get
\[
\alpha \ge \frac{{H_{i} (H_{i} - 2)}}{{H_{i} (H_{i} - 2) + 1}}.
\]
For increasing the dominance of the cubic function as much as possible, we 
use
\begin{equation}
\label{eq20}
\alpha  = \frac{{H_i (H_i  - 2)}}{{H_i (H_i  - 2) + 1}}.
\end{equation}

The proof for the case of
\[
0 < H_{i} < 1/2
\]
can be done easily as follows. With replacements of
\begin{equation}
\label{eq21}
\left\{ {{\begin{array}{*{20}c}
 {H_{i} \to 1/L_{i} ,} \hfill \\
 {k \to l - 1,} \hfill \\
\end{array} }} \right.
\end{equation}
the inequality (\ref{eq14}) is rewritten as
\[
\alpha L_{i} ^{2} + (1 - \alpha )[(1 - L_{i} )l + L_{i} ]^{3}[3(L_{i} - 1)l 
+ 2 - L_{i} ] \ge 0.
\]
This inequality corresponds to that of (\ref{eq14}) only with 
replacements of $L_{i} \to H_{i} $ and $l \to k$. Thus, the solution of this 
inequality is
\begin{equation}
\label{eq22}
\alpha \ge \frac{{L_{i} (L_{i} - 2)}}{{L_{i} (L_{i} - 2) + 1}},\quad 
{\rm for}\;L_{i} > 2\;{\rm and}\;l \in [0,1].
\end{equation}
From eq.~(\ref{eq21}), it is known that the result (\ref{eq22}) is for 
$0 < H_{i} \; < 1/2$ and $k \in [0,1]$. From the same reason as mentioned in 
leading eq.~(\ref{eq20}), we use
\begin{equation}
\label{eq23}
\alpha = \frac{{L_{i} (L_{i} - 2)}}{{L_{i} (L_{i} - 2) + 1}}.
\end{equation}

The above results on the determination of $\alpha $ is summarized as
\begin{equation}
\label{eq24}
\alpha  = \frac{{M_i (M_i  - 2)}}{{M_i (M_i  - 2) + 1}},
\end{equation}
where
\begin{eqnarray}
\label{eq25}
M_i &=& \max [2,\;\max (H_i ,\;L_i )] \nonumber \\
    &=& \max [2,\;\max (\frac{{Q_i }}{{P_i }},\;\frac{{P_i }}{{Q_i }})].
\end{eqnarray}
In the case of $1/2 \le H_{i} \le 2$ and $H_{i} \le 0$, $M_{i}$ and 
$\alpha $ are determined with eqs.~(\ref{eq25}) and (\ref{eq24}) as
\begin{eqnarray}
M_{i} = 2,\nonumber \\
\alpha  = 0,\nonumber
\end{eqnarray}
and the interpolation function is reduced to pure cubic. For $H_{i} \ge 2$, 
eq.~(\ref{eq24}) is reduced to eq.~(\ref{eq20}) and, for $L_{i} \ge 2$, 
it is reduced to eq.~(\ref{eq23}).

\subsection{Summary of the formula}
With a little rearrangement, the formula derived in the last subsection is 
summarized as
\[
F(k) = f_{i} + d_{i} hk + (G1_{i} + G2_{i} )k^{2},
\]
\[
\frac{{\partial F(k)}}{{\partial x}} = d_{i} + [G1_{i} \frac{{Q_{i} + D_{i} 
}}{{D_{i} }} + 2G2_{i} + (1 - \alpha )(Q_{i} - D_{i} )]\frac{{k}}{{h}},
\]
where
\begin{eqnarray}
G1_i &=& \alpha P_i ^2 /D_i,\nonumber \\
G2_{i} &=& (1 - \alpha )(2P_{i} - D_{i} ),\nonumber \\
D_{i} &=& Q_{i} + (P_{i} - Q_{i} )k,\nonumber \\
\alpha &=& \frac{{M_{i} (M_{i} - 2)}}{{M_{i} (M_{i} - 2) + 1}},\nonumber \\
M_{i} &=& max[2,\;max(\frac{{Q_{i} }}{{P_{i} }},\;\frac{{P_{i} }}{{Q_{i} }})]
\nonumber
\end{eqnarray}
with
\begin{eqnarray}
P_i &=& (S_i  - d_i )h,\nonumber \\
Q_{i} &=& (d_{i + 1} - S_{i} )h,\nonumber \\
S_{i} &=& (f_{i + 1} - f_{i} )/h.\nonumber
\end{eqnarray}
For the case of $u_i > 0$, we need replacements of
\[
i + 1 \to i - 1,\nonumber \quad h \to - h.
\]

\vspace{0.8cm}

\begin{figure}
\begin{center}
\epsfig{file=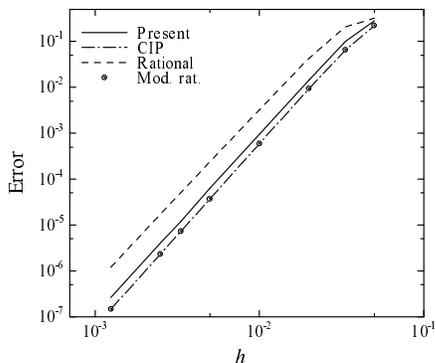,width=6cm}
\end{center}
\caption{Numerical error as a function of the grid width $h$ in the 
results with the present, the CIP, the conventional rational and the 
modified rational (with the additional switching technique) methods.}
\end{figure}

\section{NUMERICAL EXPERIMENTS}
For demonstrating the validity of the previous discussions and the accuracy 
of the present scheme, we show some numerical experiments in this section.

The first example is the linear propagation of a sinusoidal wave. The 
initial condition is
\[
f(i,0) = 0.5\cos (4\pi h\,i),
\]
where the grid width $h$ is set in this example as
\[
h = 1/N
\]
and $N$ is the number of grid points. The velocity is set as 
$u = 1$ and is assumed as constant in space and time. The boundary 
condition at $i = 0$ and $N$ is periodic. With this example, 
we compare the accuracy of the present method with that of the three 
existing methods, i.e., the CIP, the conventional rational method and that 
with the additional switching (\ref{eq8}). The last one is called below the 
modified rational method. Here we define numerical error as
\[
{\rm Error} = \frac{{\sum\limits_{i = 1}^{N} {\left| {f_{i}^{n} - 
{\rm analytical\;value}} \right|} }}{{N}},
\]
where the superscript $n$ shows the number of time step. In Fig.~1, 
we plot the error as a function of the grid width $h$ at $n = 4000$ with 
CFL = 0.2. From this result, it is known that the accuracy of the 
methods except for the conventional rational method are very similar each 
other on this problem, while that of the present one is inferior a few among 
them.

\vspace{1cm}
\begin{figure}
\begin{center}
\epsfig{file=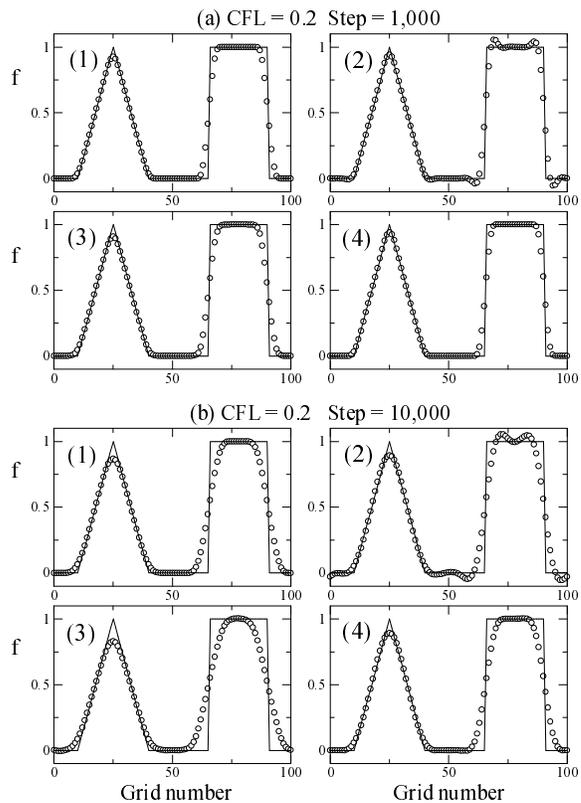,width=8cm}
\end{center}
\caption{Linear propagation of a triangular and a square waves. The figure 
shows results at $n = 1,000$ (the upper group) and $n = 10,000$ 
(the lower group) with the present (1), the CIP (2), the conventional 
rational (3) and the modified rational (4) methods.}
\end{figure}

Next we solve the linear propagation of a triangular and a square wave. The 
pulse widths of the waves are $30h$ and $26h$, respectively, 
and pulse heights are 1. In Fig.~2, we show results at $n = 1,000$ 
(a) and 10,000 (b) with CFL = 0.2. In the results with the CIP method, 
overshoots and undershoots are obviously shown around the discontinuities of 
the square wave. The results with the conventional rational method are more 
diffusive than other ones. The results with the present and the modified 
rational method are very similar while the former one is more diffusive 
quite a little and the latter one has small over- and undershoots.

With the same example of the square wave, we further discuss on the 
convexity preserving property of those methods. Here we introduce a variable 
$p$ defined as
\[
p_{i + 1/2}^{n} = \left\{ {{\begin{array}{*{20}c}
{1/2, \quad} & {{\rm if}\; f_{i + 1}^{n} - f_{i}^{n} > 0,} \\
{0, \quad} & {{\rm if}\; f_{i + 1}^{n} - f_{i}^{n} = 0,} \\
{-1/2, \quad} & {{\rm if}\; f_{i + 1}^{n} - f_{i}^{n} < 0.} \\
\end{array}} } \right.
\]
By observing this variable, we can appreciate the convexity of solution. In 
Fig.~3, we show $f$ and $p$ at $n = 150$ (a) and $10,000$ 
(b). Analytically, $p$ becomes $1/2$ and $-1/2$ only at the left and 
right discontinuities, respectively, and $p = 0$ elsewhere, namely, 
only one positive and one negative regions should exist in the spatial 
profile of $p$ if the convexity of $f$ is preserved. In the 
results with the present and the conventional rational method, the two 
regions are clearly seen while the width of the regions is expanded by the 
numerical diffusion. In the results with the CIP and the modified rational 
methods, many numbers of leaps are shown in the profile of $p$. This 
means that the conventional CIP and the modified rational methods are not 
oscillation free.

\vspace{1.2cm}
\begin{figure}
\begin{center}
\epsfig{file=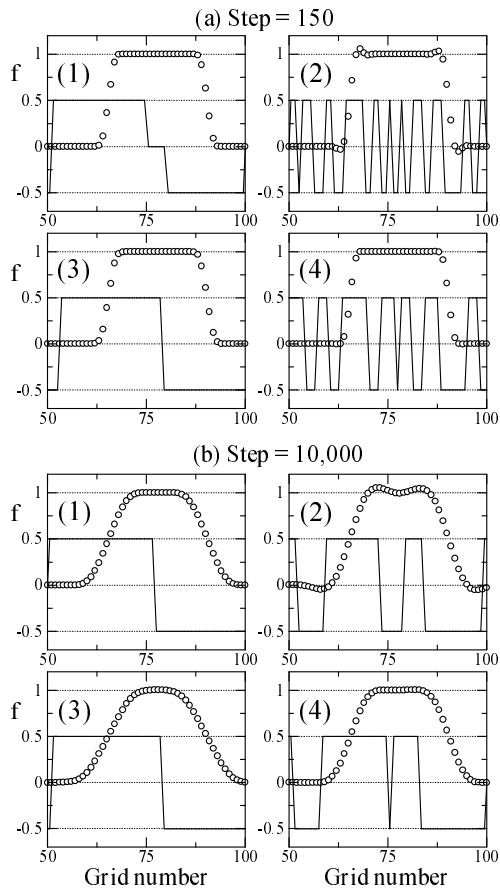,width=7cm}
\end{center}
\caption{Results of the linear propagation of a square wave at $n = 150$ 
(the upper group) and $n = 10,000$ (the lower group) with the present 
(1), the CIP (2), the conventional rational (3) and the modified rational 
(4) methods. The circles and the solid lines show $f$ and $p$, respectively.}
\end{figure}

\begin{figure}
\begin{center}
\epsfig{file=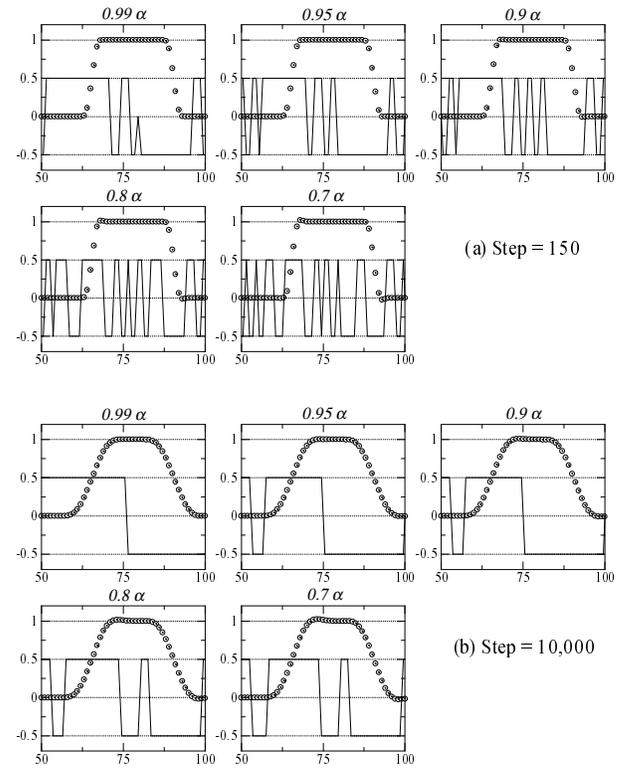,width=8.2cm}
\end{center}
\caption{Results with the reduced mixing ratio (0.99, 0.95, 0.9, 0.8 and 0.7 
times values of the optimal one) at $n = 150$ (the upper group) and 
$n = 10,000$ (the lower group).}
\end{figure}

In order to demonstrate that the determination of $\alpha $ with 
eq.~(\ref{eq24}) is the optimal, we show results of the same example with the 
present method and reduced $\alpha $. In Fig.~4, we show results 
with 0.99, 0.95, 0.9, 0.8 and 0.7 times smaller value of $\alpha $ 
than that is determined with eq.~(\ref{eq24}) and CFL = 0.2. The unphysical 
oscillation is observed even in the results with 0.99 times value and the 
oscillation becomes stronger according as $\alpha$ decreases. The 
numerical oscillation at $n = 10,000$ becomes weaker than that at 
$n = 150$ because of numerical diffusion, but still exists. The 
overshoot around the discontinuity of $f$ can also be seen in those 
results and it becomes larger as $\alpha $ decreases.

Finally we show results of an extreme example for demonstrating the robustness 
of the present method and a defect of the modified rational method. The initial 
profiles of $f$ and $u$, the latter is assumed  as constant in time, for this 
example are shown in Fig.~5. In the initial profile of $f$, there exist two 
steep gradients which are 
appropriately smooth but sufficiently sharp. In the velocity distribution, a 
similar gradient exists, and $u = 1$ in the left side of the gradient 
and $u = 0.1$ in the right side. Those mollified profiles are made as 
follows: First $f$ and $u$ are set as
\[
f_i  = \left\{ {\begin{array}{*{20}c}
{1,\quad} \hfill & {{\rm for}\;5 \le i \le 67,} \hfill  \\
{0,\quad} \hfill & {{\rm elsewhere}} \hfill  \\
\end{array}} \right.
\]
and

\begin{figure}
\begin{center}
\epsfig{file=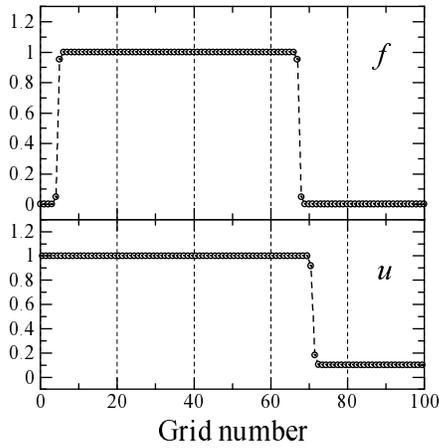,width=6cm}
\end{center}
\caption{Initial condition of the extreme example. The upper and the lower 
figures show $f$ and $u$, respectively.}
\end{figure}

\vspace{0.8cm}

\begin{figure}
\begin{center}
\epsfig{file=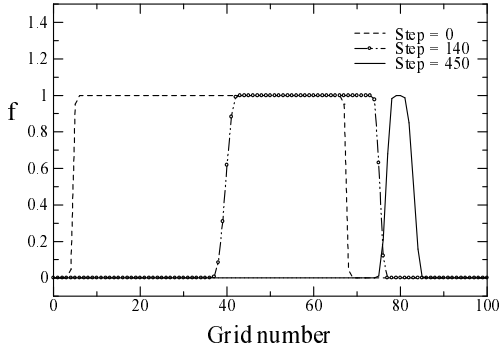,width=7cm}
\end{center}
\caption{A typical result of the extreme example by using the present method. 
The figure shows results at $n = 0$, 140 and 450 with CFL = 0.25.}
\end{figure}

\[
u_i = \left\{ {\begin{array}{*{20}c}
{1,\quad} \hfill & {{\rm for}\;i \le 71,} \hfill  \\
{0.1,\quad} \hfill & {\rm elsewhere.} \hfill  \\
\end{array}} \right.
\]
Next, those values are smoothed with a conventional 3-point smoother,
\begin{eqnarray}
g_{i}^{(m + 1)} &=& (1 - \varepsilon )g_{i}^{(m)} + \varepsilon
\frac{{g_{i + 1}^{(m)} + g_{i - 1}^{(m)} }}{{2}},\nonumber \\
&{\rm for}&\;m = 1,2, \cdots M,\nonumber
\end{eqnarray}
where $g$ is $f$ or $u$, $m$ is the iteration 
number of the smoothing, $M$ is the maximum of it and 
$\varepsilon $ is a positive constant smaller than 1.0. For 
$f$ and $u$, we use $(M, \varepsilon) = (2, 0.05)$ and 
$(2, 0.1)$, respectively. The resulting mollified values are 
used as initial values. Initial value of $d_{i}$ is set as
\[
d_i^0  = \left\{ {\begin{array}{*{20}c}
{0,\quad} \hfill & {{\rm if}\;f_i^{(M)} = 0\;{\rm or}\;f_i^{(M)} = 1,}\hfill \\
{(f_{i + 1}^{(M)} - f_{i - 1}^{(M)})/2h,\quad} \hfill & {{\rm elsewhere.}} 
\hfill \\
\end{array}} \right.
\]
In solving this example, we estimate the velocity gradient appearing in 
eq.~(\ref{eq4}) with a conventional second-order centered finite differencing 
as usually done in the CIP scheme \cite{ref1}.

\vspace{1cm}

\begin{figure}
\begin{center}
\epsfig{file=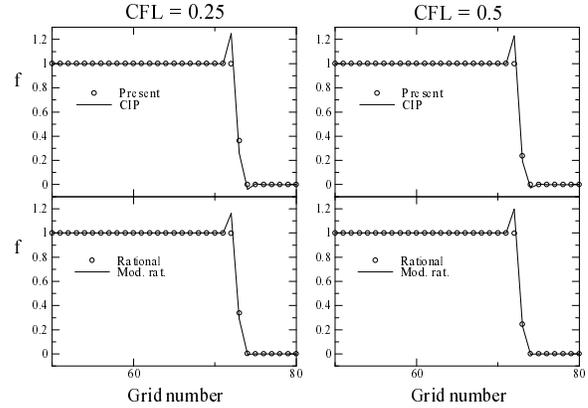,width=8cm}
\end{center}
\caption{Results when the right discontinuity just passes through the velocity 
gradient (the left group:$n = 48$ with CFL = 0.25. the right group: 
$n = 24$ with CFL = 0.5). In the results with the CIP and the modified 
rational methods, strong overshoot appears.}
\end{figure}

\vspace{1cm}

\begin{figure}
\begin{center}
\epsfig{file=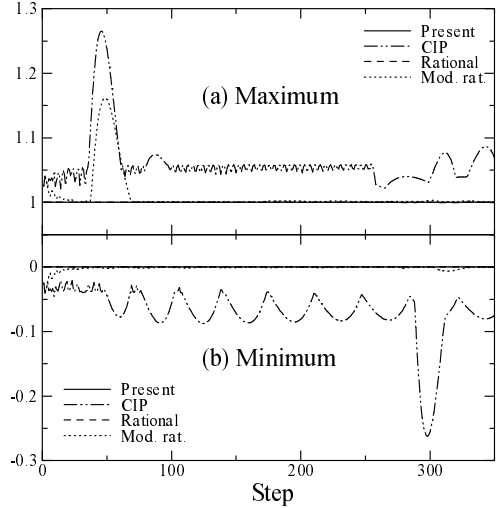,width=7cm}
\end{center}
\caption{The maximum (upper) and the minimum (lower) values in the numerical 
results of the extreme example as a function of the time step. The maximum 
and the minimum values calculated with the present and the conventional 
rational methods keep with 1 and 0, respectively, while those with the CIP 
and the modified rational method are oscillated.}
\end{figure}

In the first stage of this problem, the mollified square wave of $f$ 
propagates to right in $u = 1$. Then, after the wave covers the 
gradient of $u$, the waveform is compressed in horizontal direction 
and the pulse width is decreased. Finally, after passing through the 
velocity gradient, the width of the square becomes 1/10 of the initial one, 
i.e., $\sim 6h$, and the square propagates in $u = 0.1$ (See Fig.~6 in 
which a typical result of this problem by the present method is shown).

In Fig.~7, we show results at $n = 48$ with CFL = 0.25 and at 
$n = 24$ with CFL = 0.5. While the results with the present and the 
conventional rational methods are smooth and oscillation free, strong 
overshoot is shown in those with the CIP and the modified rational methods. 
At the time where the gradient of $f$ passes through the gradient of 
$u$, it is amplified strongly by the velocity gradient. Thus, if the 
cubic interpolation is adopted in this case, the overshoot should appear. In 
Fig.~8, we show the maximum and the minimum value of $f$ as a 
function of the number of time step. The maximums with the CIP and the 
modified rational methods have strong peak at $n \sim 50$, i.e., 
when the right gradient of $f$ 
passes through the velocity gradient and, furthermore, the minimum with the 
CIP has strong peak at $n \sim 300$, i.e., when the left gradient of 
$f$ passes through there. This result shows that the additional 
switching (\ref{eq8}) raises inadequate adaptation of the cubic function. The 
present method does not have such a defect. In Fig.~9, the results at 
$n = 550$, i.e., those after the pulse has passed through the 
velocity gradient is shown. The result with the modified rational method has 
weak overshoot and that with the conventional rational method is more 
diffusive than that with the present one.

The above results prove that the present method has higher accuracy than 
that of the conventional rational method and is a convexity-preserving 
method.

\vspace{1cm}

\begin{figure}
\begin{center}
\epsfig{file=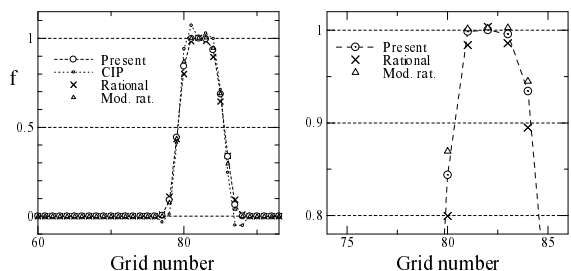,width=8cm}
\end{center}
\caption{Results after the square have passed through the velocity gradient 
($n = 550$ with CFL = 0.25). The right figure shows an enlarged one 
of a part of the left one. The result with the CIP is not shown in the right 
figure. Theoretically the height and the width of the projection should be 1 
and about $6h$, respectively.}
\end{figure}

\section{CONCLUSION}

In this paper we proposed a hybrid semi-Lagrangian method with a cubic and 
a rational interpolation functions. The optimal ratio for mixing those 
functions was led theoretically. The present method has higher accuracy than 
the conventional rational method and is oscillation free. The numerical 
experiments curried out in the last section demonstrate the validity of the 
theoretical discussion in Sec.~3 and the accuracy of the present method. 
However, the results show a limit of the present method as well. Higher 
accuracy than that of the given results shown, for example, in Figs.~2 and 3 
may no longer be expected when one use the cubic and the rational functions 
under the convexity preserving condition. Extension to a higher-order method 
may be needed for achieving higher resolution.

The present method may be extended to multidimensions and a conservative 
method by employing the approaches discussed by Aoki \cite{ref6,ref7} and 
Tanaka et al \cite{ref17,ref18}, respectively.

\end{document}